\documentclass[11pt]{article}
\usepackage{graphicx}

\newcommand{\be}{\begin{equation}}
\newcommand{\bel}[1]{\begin{equation}\label{#1}}
\newcommand{\ee}{\end{equation}}

\newcommand{\R}{\mathcal{R}}
\newcommand{\A}{\mathcal{A}}
\newcommand{\sigx}{\sigma_{\rm x}}
\newcommand{\sigy}{\sigma_{\rm y}}
\newcommand{\Mh}{\widehat{M}}
\newcommand{\eps}{\varepsilon}
\newcommand{\tr}{\mbox{tr}}
\renewcommand{\d}{\mbox{d}}
\newcommand{\px}{p_{\rm x}}
\newcommand{\py}{p_{\rm y}}
\newcommand{\Lbar}{\overline{L}}

\newcommand{\myfig}[3]{
\begin{figure}[htbp]
\centering
\resizebox{#3in}{!}{\rotatebox{0}{\includegraphics{#1.ps}}}
\caption{#2 \label{fig:#1}}
\end{figure}}

\begin{document}
\title{{\sc Bayesian alignment using hierarchical models,
with applications in protein bioinformatics}}
\author{
Peter J. Green\thanks
{School of Mathematics, University of Bristol, Bristol BS8 1TW, UK.
\newline \hspace*{5mm} Email: {\tt P.J.Green@bristol.ac.uk}.}\\
University of Bristol.\\
\and
Kanti Mardia\thanks
{School of Mathematics,
University of Leeds, Leeds, LS2 9JT, UK.
\newline \hspace*{5mm} Email: {\tt k.v.mardia@leeds.ac.uk}}\\
University of Leeds.
}
\date{Revision: 1 July 2005}
\maketitle

\begin{abstract}
An important problem in shape analysis is to match 
configurations of points in space
filtering out some geometrical transformation. In this
paper we introduce hierarchical models for such tasks, 
in which the points in the configurations are either 
unlabelled, or have at most a partial labelling constraining
the matching, and in which
some points may only appear in one of the configurations.
We derive procedures for simultaneous inference about
the matching and the transformation, using a
Bayesian approach. Our model is based on a Poisson process for
hidden true point locations; this leads to considerable 
mathematical simplification and efficiency of implementation.
We find a novel use for classic distributions from
directional statistics in a conditionally conjugate specification
for the case where the geometrical transformation includes 
an unknown rotation.
Throughout, we focus on the case of affine or rigid motion 
transformations. 
Under a broad parametric family of loss functions, an optimal Bayesian point
estimate of the matching matrix can be constructed,
that depends only on a single parameter of the family.

Our methods are illustrated by two applications from bioinformatics. 
The first problem is of matching
protein gels in 2 dimensions, and the second consists of
aligning active sites of proteins in 3 dimensions. In the latter case,
we also use information related to the grouping of the amino acids.
We discuss some open problems and suggest directions for future work.

\hspace{5mm}

\noindent {\small {\em Some key words: 
bioinformatics, 
Markov chain Monte Carlo, 
matching,
Poisson process, 
protein gels,
protein structure, 
shape analysis, 
von Mises--Fisher distribution.
}
}

\end{abstract}

\section{Introduction}
\label{sec:intro}
Various new challenging problems in shape matching have been 
appearing from different scientific areas including Bioinformatics 
and Image Analysis. 
In a class of problems in Shape Analysis,
one assumes that the points in two or more configurations are
labelled and these configurations are to be matched after filtering
out some transformation.  Usually the transformation is a rigid
transformation or similarity transformation. Several new
problems are appearing where the points of configuration are either
not labelled or the labelling is ambiguous, and in which 
some points do not appear in each of the configurations. 
An example of ambiguous labelling arises in understanding
the secondary structure of proteins, where we are
given not only the 3-dimensional molecular configuration but also the
type of molecules (amino acids) at each point. A
generic problem is to match such two configurations, where the matching has
to be invariant under some transformation group.  Descriptions of such
problems can be found in the review article by
Mardia, Taylor and Westhead (2003).

We now describe two datasets related to protein structure.  One is 
of 2-dimensional gel data where each point is a protein itself and the
transformation group is affine. In this case we have a partial
matching identified already by experts, that we can use to assess
our procedures.  
 In the second example we have a 3-dimensional configuration of
two active sites of two proteins which has also additional chemical
information. Here the underlying transformation to be filtered out is
rigid motion. In this protein structure problem, one of the main aims
is to take a query active site and find matches to a given database,
in some ranking order.  The matches will give some idea of functions of
the unknown proteins, leading to the design of new enzymes for example.

There are other related examples from Image Analysis such as matching buildings
when one has multiple 2-dimensional views of 3-dimensional objects
(see, for example, Cross and Hancock, 1998). The problem here requires
filtering out the projective transformations before matching. Other
examples involve matching outlines or
surfaces (see, for example, Chui and Rangarajan, 2000, and Pedersen, 2002).  Here 
there is no labelling of points involved, and we are dealing with a continuous
contour or surface rather than a finite number of points. Such problems
are not addressed in this paper.

In Section 2 we build a hierarchical Bayesian model for the point
configurations and derive inferential procedure for its parameters.
In particular, modelling hidden point locations as a Poisson process
leads to a considerable simplification. We discuss in particular
the problem when only a linear or affine transformation has to be
filtered out. In Section 3 we discuss prior specifications, and provide an
implementation of the resulting
methodology by means of Markov chain Monte Carlo (MCMC) samplers.
Under a broad parametric family of loss functions, an optimal Bayesian point
estimate of the matching matrix can be constructed,
which turns out to depend on a single parameter of the family.
We also discuss a modification to the likelihood in our model
to make use of partial label (`colour') information at the points.
Finally here there is a note on the possibilities
for an alternative computational
approach using the EM algorithm.
Section 4 describes application of our methods to the two examples from
Bioinformatics mentioned above: matching Protein gels in 2 dimensions and 
aligning active sites of Proteins in 3 dimensions. The paper concludes
with a Discussion of some open problems and future directions, and comparisons with
other methods. 

The principal innovations in our approach are (a) the fully model-based
approach to alignment, (b) the model formulation allowing integrating
out of the hidden point locations, (c) the prior specification for the
rotation matrix, and (d) the MCMC algorithm. 

\section{Hierarchical modelling of alignment and matching problems}
\label{sec:models}

We will build a hierarchical model for the observed point configurations, 
and derive inferential procedures for its parameters, including the
unknown matching between the configurations, according to the
Bayesian paradigm.

\subsection{Point process model, with geometrical transformation and random thinning}

Suppose we are given two point configurations in $d$-dimensional
space $\R^d$: $\{x_j, j=1,2,\ldots,m\}$ and $\{y_k, k=1,2,\ldots,n\}$.
The points are labelled for identification, but arbitrarily.

Both point sets are regarded as noisy observations on subsets of
a set of true locations
$\{\mu_i\}$, where we do not know the mappings from $j$ and $k$ to $i$.
There may be a geometrical transformation between the $x$-space
and the $y$-space, which may also be unknown.
The objective is to make model-based inference about these
mappings, and in particular make probability statements
about matching -- which pairs $(j,k)$ correspond to the same
true location?

The geometrical transformation between the $x$-space
and the $y$-space will be denoted $\A$; thus $y$ in $y$-space
corresponds to $x=\A y$ in $x$-space. The notation does not
imply that the transformation $\A$ is necessarily linear.
It may be a rotation or more general linear transformation,
a translation, both of these,
or some non-rigid motion. We regard the true locations
$\{\mu_i\}$ as being in $x$-space. 

The mappings between the indexing of the $\{\mu_i\}$ and that of
the data $\{x_j\}$ and $\{y_k\}$ are captured by indexing arrays
$\{\xi_j\}$ and $\{\eta_k\}$; specifically we assume that
\bel{likx}
x_j=\mu_{\xi_j}+\eps_{1j}
\ee
for $j=1,2,\ldots,m$, where $\{\eps_{1j}\}$ have
probability density $f_1$, and
\bel{liky}
\A y_k=\mu_{\eta_k}+\eps_{2k}
\ee
for $k=1,2,\ldots,n$, where $\{\eps_{2k}\}$ have
density $f_2$. Multiple matches are excluded, thus 
each hidden point $\mu_i$ is observed at most once in 
each of the $x$ and $y$ configurations; equivalently, the
$\xi_j$ are distinct, as are the $\eta_k$.
All $\{\eps_{1j}\}$ and $\{\eps_{2k}\}$
are independent of each other, and independent of the $\{\mu_i\}$.

\subsection{Formulation of Poisson process prior}
\label{sec:pp}

Suppose that the set of true locations $\{\mu_i\}$
forms a homogeneous Poisson process with rate $\lambda$
over a region $V\subset\R^d$ of volume $v$, and that there
are $N$ points realised in this region. Some of these 
give rise to both $x$ and $y$ points, some to points of
one kind and not the other, and some are not observed at all.
We suppose these four possibilities occur independently
for each realised point, with probabilities parameterised so
that with probabilities $(1-\px-\py-\rho\px\py,\px,
\py,\rho\px\py)$ we observe neither, $x$ alone,
$y$ alone, or both $x$ and $y$, respectively.
The parameter $\rho$ is a certain
measure of the tendency {\it a priori} for points to be matched:
the random thinnings leading to the observed $x$ and $y$ configurations
can be dependent, but remain independent from point to point. 

Given $N$, $m$ and $n$, there are $L$ matched pairs of points in our sample
if and only if the
numbers of these four kinds of occurrence among the $N$ points
are $(N-m-n+L,m-L,n-L,L)$. Under the assumptions above these four counts will
be independent Poisson distributed variables, with means
$(\lambda v (1-\px-\py-\rho\px\py),
\lambda v \px,\lambda v \py,\lambda v \rho\px\py)$.
The prior probability
distribution of $L$ conditional on $m$ and $n$ is therefore
proportional to
$$
\frac{e^{-\lambda v \px}(\lambda v \px)^{m-L}}{(m-L)!}\times
\frac{e^{-\lambda v \py}(\lambda v \py)^{n-L}}{(n-L)!}\times
\frac{e^{-\lambda v \rho\px\py}(\lambda v \rho\px\py)^L}{L!}
$$
so that
\bel{lprior}
p(L) \propto \frac{(\rho/\lambda v)^L}{(m-L)!(n-L)!L!}
\ee
for $L=0,1,\ldots,\min\{m,n\}$. The normalising constant here
is the reciprocal of 
$H(m,n,\rho/(\lambda v))$, where $H$ can be written
in terms of the confluent hypergeometric function
$$
H(m,n,d) = \frac{d^m}{m!(n-m)!} \: \mbox{}_1\!F_1(-m,n-m+1,-1/d),
$$
assuming without loss of generality that $n>m$;
see Abramowitz and Stegun (1970, p. 504).
Here and later,
we use the generic $p(\cdot)$ notation for distributions
and conditional distributions in our hierarchical model.

The matching of the configurations is represented by the
{\it matching matrix} $M$, where $M_{jk}$ indicates whether $x_j$ and $y_k$
are derived from the same $\mu_i$ point, or not, that is,
$$
M_{jk} =\cases {1 & if $\xi_j=\eta_k$ \cr
   0 & otherwise \cr}.
$$
Note that $\sum_{j,k} M_{jk}=L$, and that,
since multiple matches are ruled out, there is at most 
one 1 in each row and in each column of $M$: $\sum_j M_{jk}\leq 1 \forall k$,
$\sum_k M_{jk}\leq 1 \forall j$.
We assume for the moment that conditional
on $L$, $M$ is {\it a priori} uniform: there are $L! {m \choose L} {n \choose L}$
different $M$ matrices consistent with a given value of $L$, and these
are taken as equally likely. Thus
$$
p(M) = p(L)p(M|L) \propto \frac{(\rho/\lambda v)^L}{(m-L)!(n-L)!L!}
\left\{L! {m \choose L} {n \choose L}\right\}^{-1}
\propto (\rho/\lambda v)^L, 
$$
(where here and later `$\propto$' means proportional to, as functions of
the variable(s) to the left of the conditioning $|$, in this case, $M$).
Thus
\bel{pm}
p(M) = \frac{(\rho/\lambda v)^L}
{\sum_{\ell=0}^{\min\{m,n\}} \ell! {m \choose \ell} {n \choose \ell}(\rho/\lambda v)^\ell}.
\ee
Note that, because of the choice of parameterisation for the probabilities 
that hidden points are observed, this expression does not involve
$\px$ and $\py$.
\begin{figure}[htbp]
\centering
\resizebox{3.5in}{!}{\rotatebox{0}{\includegraphics{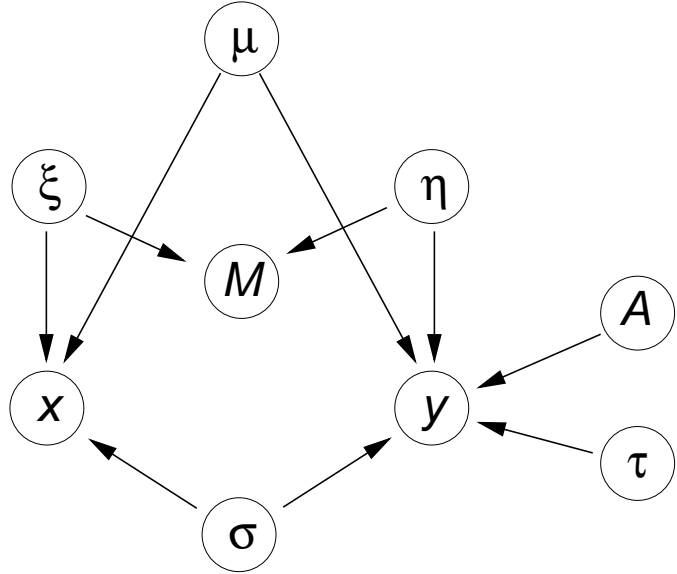}}}
\caption{Directed acyclic graph representing our model,
showing all data and parameters treated as variable. \label{fig:align}}
\end{figure}

\subsection{Likelihood of data}

We now have to specify the likelihood of the observed configurations of points,
given $M$. For simplicity, we will henceforth assume that $\A$
is an affine transformation: $\A y=Ay+\tau$. From (\ref{likx}) and
(\ref{liky}), the densities of $x_j$ and $y_k$, conditional
on $A$, $\tau$, $\{\mu_i\}$, $\{\xi_j\}$ and $\{\eta_k\}$ are
$f_1(x_j-\mu_{\xi_j})$ and $|A|f_2(Ay_k+\tau-\mu_{\eta_k})$, respectively,
$|A|$ denoting the absolute value of the determinant of $A$.

The locations $\{\mu_i\}$ of the $m-L$ points that 
generate an $x$ observation but not a $y$ observation are independently 
uniformly distributed over the region $V$, so that the likelihood contribution
of these $m-L$ observations, namely $\{x_j:M_{jk}=0\forall k\}$, is
$$
\prod_{j: M_{jk}=0\forall k} v^{-1} \int_V f_1(x_j-\mu) d\mu
$$ 
Similarly, the contributions from the unmatched $y$ observations, and from the 
matched pairs are
$$
\prod_{k: M_{jk}=0\forall j} v^{-1}  \int_V |A|f_2(Ay_k+\tau-\mu) d\mu
\quad\mbox{and}\quad
\prod_{j,k: M_{jk}=1} v^{-1} \int_V f_1(x_j-\mu) |A|f_2(Ay_k+\tau-\mu) d\mu
$$ 
respectively. These integrals all exhibit `edge effects' 
from the boundary of the region $V$, which can be neglected if
$V$ is large relative to the supports of $f_1$ and $f_2$. In this
case these three expressions approximate to
$$
v^{-(m-L)}, (|A|/v)^{n-L}, \quad\mbox{and}\quad
(|A|/v)^L \prod_{j,k: M_{jk}=1}\int_{\R^d} f_1(x_j-\mu) f_2(Ay_k+\tau-\mu) d\mu
$$
respectively. The last expression can be written
$$
(|A|/v)^L \prod_{j,k: M_{jk}=1} g(x_j-Ay_k-\tau) 
$$
where $g(z)=\int f_1(z+u)f_2(u)du$ (the density of $\eps_{1j}-\eps_{2k}$).

Combining these terms, the complete likelihood is
\bel{lik}
p(x,y|M,\A) = v^{-(m+n)} |A|^n \prod_{j,k: M_{jk}=1} g(x_j-Ay_k-\tau).
\ee
Multiplying (\ref{pm}) and (\ref{lik}), we then have
$$
p(M,x,y|\A) \propto |A|^n \prod_{j,k: M_{jk}=1} \{(\rho/\lambda)g(x_j-Ay_k-\tau)\}.
$$
Note that the constant of proportionality involves $m$, $n$, $\lambda$, $\rho$,
and $v$, but not $A$, $\tau$, any parameters in $f_1$ or $f_2$, or $M$ of course.

If we further specialise by making assumptions of spherical normality
for $f_1$ and $f_2$:
$$
x_j \sim N_d(\mu_{\xi_j},\sigx^2I)
\qquad\mbox{and}\qquad
Ay_k+\tau \sim N_d(\mu_{\eta_k},\sigy^2I),
$$
with $\sigx=\sigy=\sigma$, say, then
$$
g(z)=\frac{1}{(\sigma\surd 2)^d} \phi(z/\sigma\surd 2)
$$
where $\phi$ is the standard normal density in $\R^d$, and 
our final joint model is
\bel{post}
p(M,A,\tau,\sigma,x,y) \propto
|A|^{n}  p(A) p(\tau) p(\sigma)\prod_{j,k:M_{jk}=1} \left(
\frac{\rho\phi(\{x_j-Ay_k-\tau\}/\sigma\surd 2)}{\lambda(\sigma\surd 2)^d}\right).
\ee
Note that not only $\px$ and $\py$ but also $v$ does not appear in this
expression, principally from our choice of parameterisation, 
and that only the ratio $\rho/\lambda$ is identifiable. The directed acyclic
graph representing this joint probability model, including the
variables ($\mu$, $\xi$ and $\eta$) that we have integrated out, is displayed in
Figure \ref{fig:align}.

\section{Prior distributions and computational implementation}
\label{sec:mcmc}

We will henceforth treat $\rho$ and $\lambda$ as fixed,
and consider inference for the remaining unknowns
$M$, $\tau$, $\sigma^2$ and sometimes $A$, given
the data $\{x_j\}$ and $\{y_k\}$. Markov chain Monte Carlo
methods must be used for the computation; several introductions
and overviews of MCMC are available, for example, the primer
in Green (2001). In Section \ref{sec:em}, we discuss the relevance
and applicability of an EM algorithm for making inference with
an approximation of our model.

We suppose
that prior information about $\tau$, $\sigma^2$ and $A$ will be
at best weak, and so we concentrate on generic prior
formulations that facilitate the posterior analysis.
Prior assumptions are therefore discussed in parallel 
with MCMC implementation. Note that our formulation has some
affinity with mixture models, the matching matrix $M$ playing a similar
role to the allocation variables often used in computing
with mixtures; see, for example, Richardson and Green (1997).
As in that paper, the fully Bayesian analysis here aims at
simultaneous joint inference about both the discrete and
continuously varying unknowns, in contrast to frequentist
approaches. 

Our model has another similarity with a mixture
formulation, in that as $M$ varies, the number of
hidden points needed to generate all the observed
data also varies, and thus there seems to be a
`variable-dimension' aspect to the model.
However, here our approach of integrating out the
hidden point locations eliminates the variable-dimension
parameter, so that reversible jump MCMC is not needed.

\subsection{Priors and MCMC updating for a rotation matrix}
\label{sec:rotmat}

We are interested in alignment and matching problems in which
either $A$ is given, and treated as fixed, or in which it is
one of the objects of inference. In the latter case, we consider in this
paper only the case of rotation matrices in two and three
dimensions. We therefore focus on the full
conditional distribution for $A$, which from (\ref{post}) is
$$
p(A|M,\tau,\sigma,x,y) \propto
|A|^{n}  p(A) \prod_{j,k:M_{jk}=1}
\phi(\{x_j-Ay_k-\tau\}/\sigma\surd 2).
$$
Viewing this as a density for $A$, we are still free to choose the
dominating measure for $p(A)$, which is arbitrary: this full conditional
density is then with respect to the same measure.

Let us restrict attention to {\it rotations}: orthogonal matrices $A$,
(those with $A^{-1}$ = $A^T$) with positive determinant,
so that $|A|=1$.
Expanding the expression above, we then find
$$
p(A|M,\tau,\sigma,x,y) \propto
  p(A) \exp\left(\sum_{j,k:M_{jk}=1} 
-0.5(||x_j-Ay_k-\tau||/\sigma\surd 2)^2
\right)
$$
$$
\propto
p(A) \exp \left(\tr\left\{ (1/2\sigma^2)\sum_{j,k:M_{jk}=1} y_k(x_j-\tau)^TA\right\} \right).
$$

Note a remarkable opportunity for (conditional) conjugacy -- if $p(A)$
has the form $p(A)\propto \exp(\tr(F_0^TA))$ for some matrix $F_0$, then
the posterior has the same form with $F_0$ replaced by
$$
F=F_0+(1/2\sigma^2)\sum_{j,k:M_{jk}=1} (x_j-\tau)y_k^T.
$$
This form of $p(A)$ is known as the matrix Fisher distribution (Downs, 1972; Mardia and Jupp, 2000, p. 289).
To the best of our knowledge, this unique role of the matrix Fisher distribution 
(or in the two-dimensional case, the von Mises distribution) as the 
prior distribution for a rotation conjugate to spherical
Gaussian error distributions has not previously been noted.
(Although Mardia and El-Atoum (1976) have identified the von Mises--Fisher
distribution as the conjugate prior for the mean direction).
This may have relevance in models for other situations, 
including the simpler case where there is no uncertainty in the
matching. The conjugacy is presumably related to the
interpretation of the matrix Fisher distribution
as a conditional multivariate Gaussian (see Mardia and Jupp, 
2000, p.289).

\subsubsection*{Two-dimensional case}
Now consider the two-dimensional case, $d=2$. An arbitrary
rotation matrix $A$ can be written
$$
A=\left(
\begin{array}{rr}
\cos \theta &  -\sin \theta \\
\sin \theta&  \cos \theta \\
\end{array}
\right)
$$
and the natural dominating measure for $\theta$ is 
Lebesgue on $(0,2\pi)$. Then
a uniformly distributed choice of $A$ corresponds to $p(A)\propto 1$.
More generally, the von Mises distribution for $\theta$
$$
p(\theta) \propto \exp(\kappa\cos(\theta-\nu))=\exp(\kappa\cos\nu\cos\theta+\kappa\sin\nu\sin\theta)
$$
can indeed be expressed as $p(A)\propto \exp(\tr(F_0^TA))$, where
a (non-unique) choice for $F_0$ is
$$
F_0=\kappa/2\left(
\begin{array}{rr}
\cos \nu &  -\sin \nu \\
\sin \nu&  \cos \nu \\
\end{array}
\right).
$$

Thus the full conditional distribution for $\theta$ is of the same
von Mises form, with
$\kappa\cos\nu$ updated to $(\kappa\cos\nu+S_{11}+S_{22})$, and
$\kappa\sin\nu$ to $(\kappa\sin\nu-S_{12}+S_{21})$,
where $S$ is the $2\times 2$ matrix $(1/2\sigma^2)\sum_{j,k:M_{jk}=1} (x_j-\tau)y_k^T$.

It is therefore trivial to implement a Gibbs sampler
move to allow inference about $A$, assuming a von Mises 
prior distribution on the rotation angle $\theta$
(including the uniform case, $\kappa=0$). We can use the
Best/Fisher algorithm, an efficient rejection method
(see Mardia and Jupp, 2000, p.43),
to sample from the full conditional for $\theta$.

\subsubsection*{Three-dimensional case}

In the three-dimensional case, we can represent $A$ as the
product of elementary rotations
\bel{geneul}
A=A_{12}(\theta_{12})A_{13}(\theta_{13})A_{23}(\theta_{23})
\ee
as in Raffenetti and Ruedenberg (1970), and Khatri and Mardia (1977).
Here, for $i<j$, $A_{ij}(\theta_{ij})$ is the matrix with 
$m_{ii}=m_{jj}=\cos\theta_{ij}$, $-m_{ij}=m_{ji}=\sin\theta_{ij}$,
$m_{rr}=1$ for $r\neq i,j$ and other entries 0.
We can then update each of the generalised Euler angles
$\theta_{ij}$ in turn, conditioning on the other two angles and the
other variables ($M,\tau,\sigma,x,y$) entering the expression for
$F$. 

The joint full conditional density of the Euler angles is
$$
\propto \exp[\tr\{F^TA\}] \cos\theta_{13}
$$
for $\theta_{12},\theta_{23}\in(-\pi,\pi)$ and $\theta_{13}\in(-\pi/2,\pi/2)$.
The cosine term arises since the natural dominating measure,
corresponding to uniform distribution of rotation,
has volume element $\cos\theta_{13} \:\d\theta_{12}\:\d\theta_{13}\:\d\theta_{23}$
in these coordinates.

Substituting the representation (\ref{geneul}), and simplifying, we
find that the trace can be written variously as
$\tr\{F^TA\}=a_{12}\cos\theta_{12}+b_{12}\sin\theta_{12}+c_{12}
=a_{13}\cos\theta_{13}+b_{13}\sin\theta_{13}+c_{13}
=a_{23}\cos\theta_{23}+b_{23}\sin\theta_{23}+c_{23}$
where
\begin{eqnarray*}
a_{12}& = &(F_{22}-\sin\theta_{13}F_{13})\cos\theta_{23}
+(-F_{23}-\sin\theta_{13}F_{12})\sin\theta_{23}
+\cos\theta_{13}F_{11}
\\
b_{12}& = &(-\sin\theta_{13}F_{23}-F_{12})\cos\theta_{23}
+(F_{13}-\sin\theta_{13}F_{22})\sin\theta_{23}
+\cos\theta_{13}F_{21}
\\
a_{13}& = &\sin\theta_{12}F_{21}+\cos\theta_{12}F_{11}+\sin\theta_{23}F_{{32}}+\cos\theta_{23}F_{33}
\\
b_{13}& = &
(-\sin\theta_{23}F_{12}-\cos\theta_{23}F_{13})\cos\theta_{12}
+(-\sin\theta_{23}F_{22}-\cos\theta_{23}F_{23})\sin\theta_{12}+F_{31}
\\
a_{23}& = &(F_{22}-\sin\theta_{13}F_{13})\cos\theta_{12}+(-\sin\theta_{13}F_{23}-F_{12})\sin\theta_{12}
+\cos\theta_{13}F_{33}
\\
b_{23}& = &(-F_{23}-\sin\theta_{13}F_{12})\cos\theta_{12}+(F_{13}-\sin\theta_{13}F_{22})\sin\theta_{12}
+\cos\theta_{13}F_{32}
\end{eqnarray*}
and the $c_{ij}$ can be ignored, combined into the normalising constants.
Thus the full conditionals for $\theta_{12}$ and $\theta_{23}$ are von Mises distributions,
and so these two variables can be updated by Gibbs sampling.
That of $\theta_{13}$ is 
proportional to
$$
\exp[a_{13}\cos\theta_{13}+b_{13}\sin\theta_{13}]\cos\theta_{13}
$$
and we use a random walk Metropolis update for this variable,
with a perturbation uniformly distributed on $[-0.1,0.1]$.
The latter distribution has been studied in Mardia and 
Gadsden (1977) but with no discussion on how to 
simulate from it.

\subsection{Priors and updating for other parameters}
\label{sec:prior}

We make the standard normal/inverse gamma assumptions:
$$
\tau \sim N_d(\mu_\tau,\sigma_\tau^2I)
\qquad\mbox{and}\qquad
\sigma^{-2} \sim \Gamma(\alpha,\beta).
$$

Under the assumptions of (\ref{post}), there
is conjugacy for $\tau$ and $\sigma$, and we have explicit
full conditionals:
$$
\tau|M,A,\sigma,x,y \sim N_d\left(
\frac{\mu_\tau/\sigma_\tau^2+\sum_{j,k:M_{jk}=1} (x_j-Ay_k)/2\sigma^2}
{1/\sigma_\tau^2+L/2\sigma^2}
,
\frac{1}{1/\sigma_\tau^2+L/2\sigma^2}I
\right)
$$
$$
\sigma^{-2}|M, A,\tau,x,y \sim \Gamma\left(\alpha+(d/2)L,
\beta+(1/4)\sum_{j,k:M_{jk}=1} ||x_j-Ay_k-\tau||^2\right),
$$
and so it is trivial to implement Gibbs sampler updates for these parameters.

\subsection{Updating $M$}
\label{sec:updatem}
The matching matrix $M$ is updated in detailed balance using Metropolis-Hastings
moves that only propose changes to a few entries: the number of matches
$L=\sum_{j,k}M_{jk}$ can only increase or decrease by 1 at a time, or 
stay the same. The possible changes are 
\begin{enumerate}
\item[(a)] adding a match: changing one entry $M_{jk}$ from 0 to 1
\item[(b)] deleting a match: changing one entry $M_{jk}$ from 1 to 0
\item[(c)] switching a match: simultaneously changing one entry from 0 to 1, and another {\it in the same row or column} from 1 to 0.
\end{enumerate}

The proposal proceeds as follows: first a uniform random choice is made
from all the $m+n$ data points $x_1,x_2,\ldots,x_m,y_1,y_2,\ldots,y_n$.
Suppose without loss of generality, by the symmetry of the set-up, that
an $x$ is chosen, say $x_j$. There are two possibilities: either $x_j$ is
currently matched ($\exists k$ such that $M_{jk}=1$) or not (there is no
such $k$).

If $x_j$ is matched to $y_k$, with probability $p^\star$ we
propose {\it deleting} the match, and with probability $1-p^\star$ we 
propose {\it switching}
it from $y_k$ to $y_{k'}$, where $k'$ is drawn uniformly at random from
the currently unmatched $y$ points. On the other hand, if $x_j$
is not currently matched, we propose {\it adding} a match between
$x_j$ and a $y_{k}$, where again $k$ is drawn uniformly at random from
the currently unmatched $y$ points.

The acceptance probabilities for these three possibilities are easily
derived from the expression (\ref{post}) for the joint distribution, since in each case the proposed new matching matrix $M'$ is only slightly perturbed from $M$,
so that the ratio $p(M',\tau,\sigma|x,y)/p(M,\tau,\sigma|x,y)$ has only 
a few factors. Taking into account also the proposal probabilities, 
whose ratio is $(1/n_{\rm u})\div p^\star$,
where $n_{\rm u}=\#\{k\in 1,2,\ldots,n: M_{jk}=0\forall j\}$ is the number of
unmatched $y$ points in $M$, we find
that the acceptance probability for 
adding a match $(j,k)$ is
\bel{mhadd}
\min\left\{1,
\frac
{\rho\phi(\{x_j-Ay_k-\tau\}/\sigma\surd{2})p^\star n_{\rm u}}
{\lambda(\sigma\surd{2})^d}
\right\}.
\ee
Similarly, the acceptance probability for switching the match of $x_j$ from
$y_k$ to $y_{k'}$ is
\bel{mhswitch}
\min\left\{1,
\frac{\phi(\{x_j-Ay_{k'}-\tau\}/\sigma\surd{2})}
{\phi(\{x_j-Ay_k-\tau\}/\sigma\surd{2})}\right\}
\ee
and for 
deleting the match $(j,k)$ it is
$$
\min\left\{1,
\frac{\lambda(\sigma\surd{2})^d}
{\rho\phi(\{x_j-Ay_k-\tau\}/\sigma\surd{2})p^\star n_{\rm u}'}\right\},
$$
where $n_{\rm u}'=\#\{k\in 1,2,\ldots,n: M_{jk}'=0\forall j\}=n_{\rm u}+1$.
Along with just one of each of the other updates,
we typically make several moves updating $M$ per sweep,
since the changes effected are so modest.

\subsection{Loss functions}

The output from the MCMC sampler derived above, once equilibrated,
is a sample from the posterior distribution determined by (\ref{post}). 
As always with sample-based computation, this 
provides an extremely flexible basis for reporting
aspects of the full joint posterior that are of interest.

The matching matrix $M$ will often be of particular inferential
interest, and for some purposes a point estimate is desirable; in this
section we discuss how to obtain a Bayesian point estimate of the 
matching matrix $M$. 

The most easily understood estimator of $M$ would
be its posterior mode, the {\it maximum a posteriori} (MAP)
estimator. However, there are difficulties here. First, the notion is
itself ambiguous -- the unknown `parameter' in our model
consists of the matching matrix $M$, and some real parameters.
`MAP' might refer to the $M$ component of the overall maximum,
or the mode of the marginal posterior for $M$ alone. Secondly,
the posterior is multi-modal, and different modes may have
different `widths', appropriately measured. So there is no
intrinsic attraction to the MAP estimate. We should return
to basic principles.

By standard theory, this requires specification of
a loss function, $L(M,\widehat{M})$, giving the cost
incurred in declaring the matching matrix to be
$\widehat{M}$ when it is in fact $M$. The optimal estimate
given data $(x,y)$ is the matching matrix $\widehat{M}$ that
minimises the posterior expected loss
$$
E[L(M,\widehat{M})|x,y],
$$
the expectation over $M$ being taken with respect to the posterior
determined by (\ref{post}). In this language, the MAP estimator
is optimal for the `zero--one' loss function under which a fixed
total cost is paid if there is a single error in any value
$M_{jk}$; this is logically unappealing, and a further
argument against using MAP.

We consider instead loss functions $L(M,\widehat{M})$ that penalise different
kinds of error and do so cumulatively.
The simplest of these are additive over pairs $(j,k)$. Suppose that the loss when
$M_{jk}=a$ and $\widehat{M}_{jk}=b$, for $a,b=0,1$ 
is $\ell_{ab}$; for example, $\ell_{01}$ is the loss associated
with declaring a match between $x_j$ and $y_k$ when there is really none,
that is, a `false positive'. Then it is readily shown that
$$
E[L(M,\widehat{M})|x,y] = -(\ell_{10}+\ell_{01}-\ell_{11}-\ell_{00})
\sum_{j,k:\Mh_{jk}=1} (p_{jk}-K)
$$
where  
$$
K=(\ell_{01}-\ell_{00})/
(\ell_{10}+\ell_{01}-\ell_{11}-\ell_{00}),
$$
and 
$p_{jk}=p(M_{jk}=1|x,y)$ is the posterior probability that $(j,k)$
is a match, which is estimated from an MCMC run by the 
empirical frequency of this match.
Thus, provided that $\ell_{10}+\ell_{01}-\ell_{11}-\ell_{00}>0$ and
$\ell_{01}-\ell_{00}>0$, as is natural, the optimal estimate is that maximising
the sum of marginal posterior probabilities of the declared matches
$\sum_{j,k:\Mh_{jk}=1} p_{jk}$, penalised by a multiple $K$ times the number of matches.
The optimal match therefore depends on the four loss function parameters
only through the cost ratio $K$. If false positive and false negative matches 
are equally undesirable, one can simply choose $K=0.5$.

Computation of the optimal match $\Mh$ would be trivial but for the constraint that
there can be at most one positive entry in each row and column of the array.
For modest-sized problems, the optimal match can be found by informal
heuristic methods. These may not even be necessary, especially if $K$ is not too small.
In particular, it is immediate that if the set of all $(j,k)$ pairs
for which $p_{jk}>K$
includes no duplicated $j$ or $k$ values, the optimal $\Mh$ consists of precisely
these pairs.

We could also consider loss functions that penalise mismatches differently from
the sum of the losses of the individual errors. For example, declaring
$(j,k)$ to be a match when it should be $(j,k')$ might deserve a relative
loss greater or lesser
than $(\ell_{10}+\ell_{01}-\ell_{11}-\ell_{00})$, depending
on context. Such loss functions
could be handled in a broadly similar way, but this is left for future work.

\subsection{Using partial labelling information}
\label{modlik}
When the points in each configuration are `coloured', with the
interpretation that like-coloured points are more likely to be matched
than unlike-coloured ones, it is appropriate to use a modified likelihood
that allows us to exploit such information. Let the colours for the $x$ and $y$
points be $\{r^{\rm  x}_j,j=1,2,\ldots,m\}$ and 
$\{r^{\rm y}_k,k=1,2,\ldots,n\}$ respectively. The hidden point
model is augmented to generate the point colours, as follows.
Independently for each hidden point, with probability
$(1-\px-\py-\rho\px\py)$ we observe neither
$x$ nor $y$ point, as before. With probabilities
$\px\pi^{\rm x}_r$ and  $\py\pi^{\rm y}_r$,
respectively, we observe only an $x$ or $y$ point, 
with colour $r$ from an appropriate finite set.
With probability 
$$
\rho\px\py\pi^{\rm x}_r\pi^{\rm y}_s
\exp\{\gamma I[r=s]+ \delta I[r \neq s]\},
$$
we observe an $x$ point coloured $r$ and a $y$ point
coloured $s$.
Our original likelihood is equivalent to the case $\gamma=\delta=0$,
where colours are independent and so carry no information about matching. 
If $\gamma$ and $\delta$ increase, then matches are more probable, 
{\it a posteriori}, and if $\gamma>\delta$, matches between like-coloured 
points are more likely than those between unlike-coloured ones.
The case $\delta\to-\infty$ allows the prohibition
of matches between unlike-coloured points, a feature that might be adapted
to other contexts such as the matching of shapes with
given landmarks.

In implementation of this modified likelihood, the
MCMC acceptance ratios in Section \ref{sec:updatem}
have to be modified accordingly. For example,
if $r^{\rm x}_j=r^{\rm y}_k$ and $r^{\rm x}_j\neq r^{\rm y}_{k'}$,
then (\ref{mhadd}) has to be multiplied by
$\exp(-\gamma)$ and (\ref{mhswitch}) by $\exp(\delta-\gamma)$.

Other, more complicated, colouring distributions where the log probability
can be expressed linearly in entries of $M$ can be handled 
similarly.

\subsection{Alternative approach using the EM algorithm}
\label{sec:em}

The interplay between matching (allocation) and parameter uncertainty
has something in common with mixture estimation. This might suggest considering
maximisation of the posterior by using the EM algorithm,
which could of course in principle be applied either to maximum likelihood estimation
based on (\ref{lik}) or to MAP estimation based on (\ref{post}).
For the EM formulation, the `missing data' are the matches. 

In an 
exponential family, the EM algorithm alternates between
between finding expectations of missing values 
given data, at current parameter values, and
maximising the log-posterior, with missing values
replaced by these expectations.

The `expectations of missing values' are just probabilities of
matching. These are only tractable if we were to drop the
assumption that a point can only be matched with at most one other
point -- that is, that $\sum_j M_{jk}\leq 1 \forall k$, $\sum_k M_{jk} \leq 1 \forall j$.
Making this approximation, the E-step is trivial:
the expectation of $I[M_{jk}=1]$ is 
$p_{jk}=w_{jk}/(1+w_{jk})$ where $w_{jk}$ is the $(j,k)$ factor in the joint model, i.e.
$$
w_{jk}=\{(\rho/\lambda)g_\sigma(x_j-Ay_k-\tau)\}
$$

The M-step then requires maximising (for given $p_{jk}$)
$$
\log\left[|A|^{n}  p(A) p(\tau) p(\sigma)\right]
+\sum_{j,k} p_{jk} \log \{w_{jk}(A,\tau,\sigma)\}
$$
over $A$, $\tau$, $\sigma$ -- note that here $w_{jk}$ is a function
of all three. Although for some individual parameters this seems to
be explicit, in the general case we need numerical optimisation.

In summary, EM allows us to study only certain aspects of an approximate
version of our model, and is not trivial numerically -- so we
do not pursue this approach. Obtaining the complete posterior by MCMC sampling
gives much greater freedom in inference.

\section{Applications}

\subsection{Matching protein gels}

The objective in this example is to match two electrophoretic gels
automatically, given the locations of the centres of 35
proteins on each of the two gels. The data are presented
in the supplementary information on the web.
The correspondence between pairs of proteins, one
protein from each gel, is unknown, so our aim is to match the two gels
based on these sets of unlabelled points. We suppose it
is known that the transformation between the gels is affine. 
In this case, experts have already identified 10
points; see Horgan et al (1992).
Based on these 10 matches, the linear part of the
transformation is estimated {\it a priori} to be
\bel{gelA}
A=\left(
\begin{array}{rr}
0.9731 &  0.0394 \\
-0.0231 &   0.9040 \\
\end{array}
\right).
\ee
(Dryden and Mardia, 1998, pp. 20--21, 292--296).

\begin{table}
\caption{The 20 marginally most probable matches in the analysis of the gel data.
\label{gelmatches}}
\footnotesize
\vspace*{5mm}\centering\leavevmode
\begin{tabular}{crrl}
rank & $j$ & $k$ & $p_{jk}$ \\
\hline
1 & 15 & 21 & 1 \\
2 & 19 & 19 & 1 \\
3 & 8 & 8 & 1 \\
4 & 3 & 3 & 1 \\
5 & 2 & 2 & 1 \\
6 & 31 & 30 & 0.9989 \\
7 & 6 & 6 & 0.9987 \\
8 & 4 & 4 & 0.9966 \\
9 & 5 & 5 & 0.9946 \\
10 & 10 & 10 & 0.9927 \\
11 & 24 & 23 & 0.9855 \\
12 & 7 & 7 & 0.9824 \\
13 & 32 & 31 & 0.9776 \\
14 & 1 & 1 & 0.9763 \\
15 & 9 & 9 & 0.9677 \\
16 & 26 & 32 & 0.7910 \\
17 & 12 & 13 & 0.7552 \\
18 & 21 & 33 & 0.3998 \\
19 & 26 & 27 & 0.1931 \\
20 & 35 & 35 & 0.0025 \\
\hline
\end{tabular}
\end{table}
\myfig{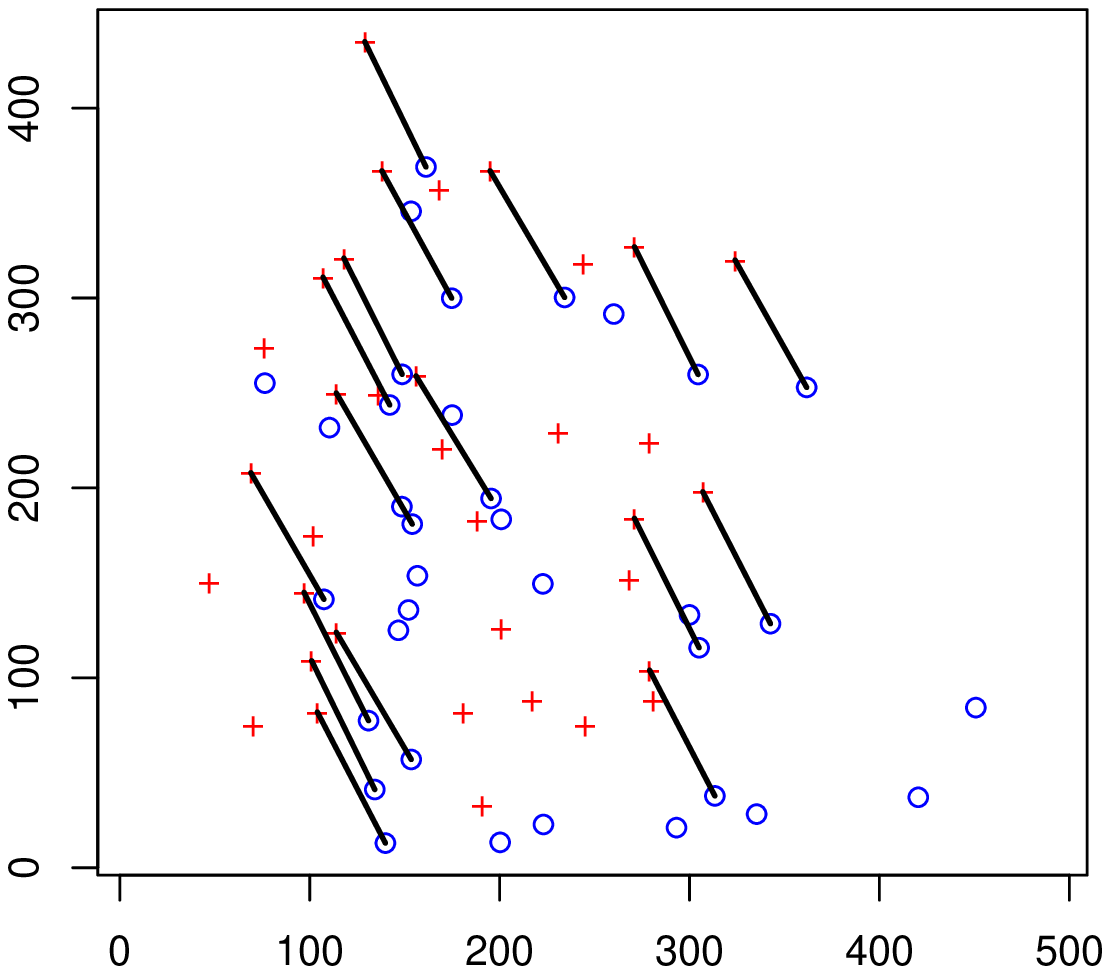}{The 17 most probable matches in the gel data, the optimal
match for any $K\in(0.3998,0.7552)$;
+ symbols signify $x$ points, o symbols the $y$ points, linearly transformed 
by premultiplication by the fixed affine transformation $A$ given in (\ref{gelA}).
The solid line for each of the 17 matches
joins the matched points, and represents the inferred translation $\tau$
plus noise.}{5}

Here, we have only
to make inference on the translation $\tau$ and the unknown matching
between certain of the proteins. 
The model (\ref{post}) will therefore be taken to apply, with $d=2$
and with $A$ held fixed at (\ref{gelA}).
The MCMC sampler described in Section \ref{sec:mcmc} was run for 
100 000 sweeps, after a burn-in period of 20 000 sweeps, considered on the
basis of an informal visual assessment of time series traces to be adequate
for convergence. 
Prior and hyperprior settings were: $\alpha=1$, $\beta=16$, $\mu_\tau=(0,0)^T$,
$\sigma_\tau=20.0$ and $\lambda/\rho=0.0001$. The sampler parameter $p^\star$ was set to 0.5.
Such a run took about 2 seconds on a 800MHz PC. Acceptance rates for the moves updating
$M$ were between 0.6\% and 2.1\%.

The posterior expectation and variance of $\tau$ were estimated to
be $(-35.950,66.685)^T$ 
(to be compared with $(-36.08,66.64)^T$ obtained by Dryden and Mardia (1998))
and
$$
\left(
\begin{array}{rr}
 0.5776 & -0.0227 \\
-0.0227 &  0.6345 \\
\end{array}
\right).
$$
The posterior mean and variance of $\sigma$ are 2.050 and 0.1192.

The 20 most probable matches between $x$ and $y$ points are listed in
Table \ref{gelmatches};
note that there is no duplication in their indices until the 19th
match: $j=26$ also appears in the 16th match (recall that there is a 
simple rule for identifying the optimal $\widehat{M}$
if there are no duplicates among the matches with $p_{jk}$ above the threshold $K$). 
We can conclude that for all values of
$K=(\ell_{01}-\ell_{00})/
(\ell_{10}+\ell_{01}-\ell_{11}-\ell_{00})$
from 1 down to 0.1112, the optimal Bayesian matching is given
by an appropriate subset of Table \ref{gelmatches}, reading down from the top.
For example if this cost ratio is 0.8 we take the first 15
rows of the table, while if the ratio is 0.6 or 0.4 we include the 16th and 17th rows
as well. The 17 most probable matches are displayed graphically
in Figure \ref{fig:gelconf}.

It will be noted that all of the expert-identified matches, points 1 to
10 in each set, are declared to be matches with high probability in the
Bayesian analysis. We also repeated the analysis with these 10 pairs
held fixed. The next 9 most probable matches, together with these 10,
are identical to those in the first 19 lines of Table \ref{gelmatches}, and
the posterior probabilities differ by less than 0.037 in all 19 cases.


\subsection{Aligning proteins in three dimensions}
\label{sec:3deg}
We now apply the matching method to a problem in three dimensional
structural biology, previously considered by Gold et al (2002).
The problem consists of finding the matches for two Active sites 1 and 2 corresponding to two Proteins A and B respectively. The corresponding coordinates $x$ and $y$ of these sites
are presented in the supplementary information; these coordinates are the centres of gravity
of the amino acids of the two sites. Here $m= 40$ and $n=63$. The biological details of the two proteins are as follows. Protein 1 is the human protein `17--beta hydroxysteroid dehydrogenase' and is involved in the synthesis of oestrogens.  This protein binds the ligands (molecules comparatively smaller than proteins) oestradiol and NADP. Protein 2 is the mouse protein `carbonyl reductase' and is involved in metabolism of carbonyl compounds.  This protein binds the ligands 2--Propanol and NADP. The common element between these two sets of ligands is NADP. From chemical properties of the sites, the relevant matching should be invariant under rigid transformation. 


\label{sec:3dalign}

\myfig{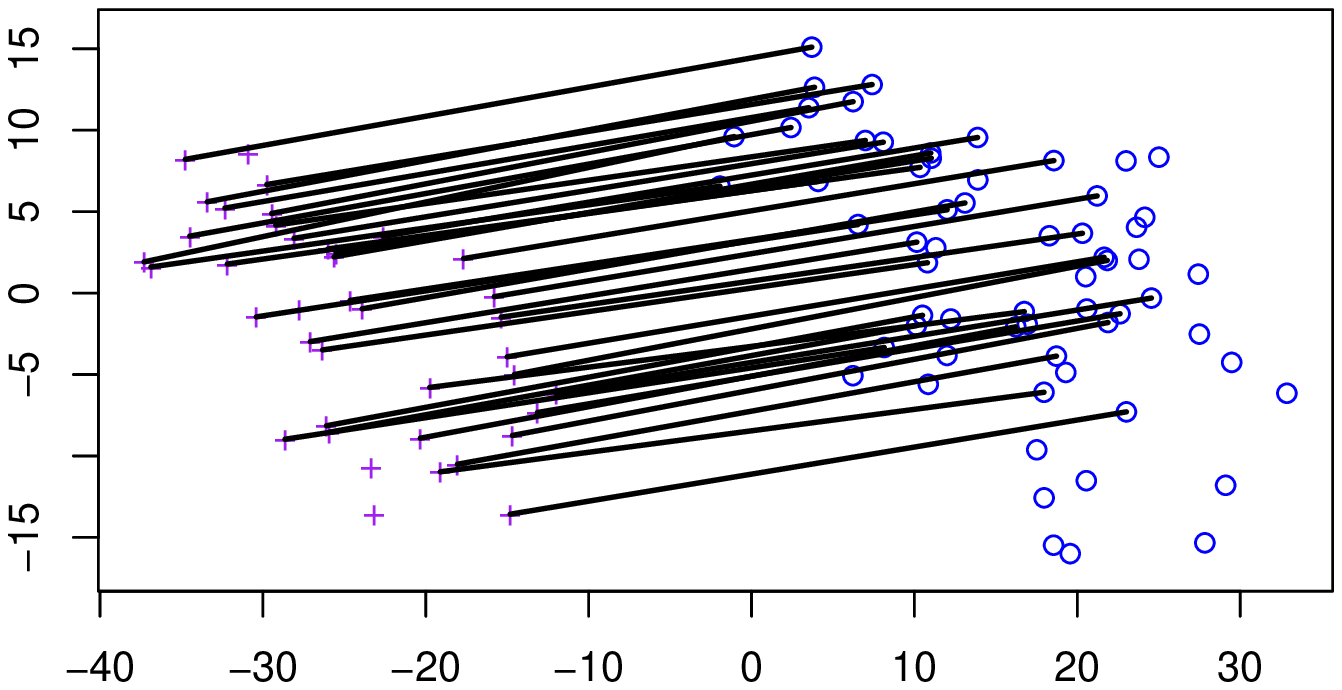}{The optimal alignment (36 matches)
when $K=0.5$ for the protein alignment analysis data, without using 
colouring information; + symbols signify $x$ points, 
o symbols the $y$ points, rotated according to the 
inferred $\widehat{A}$ matrix given by (\ref{meana1}). The entire joint configuration
has been rotated to its first two principal axes.
Solid lines represent the 36 marginally most probable matches, and 
indicate the inferred translation $\tau$ plus noise.}{5}

\myfig{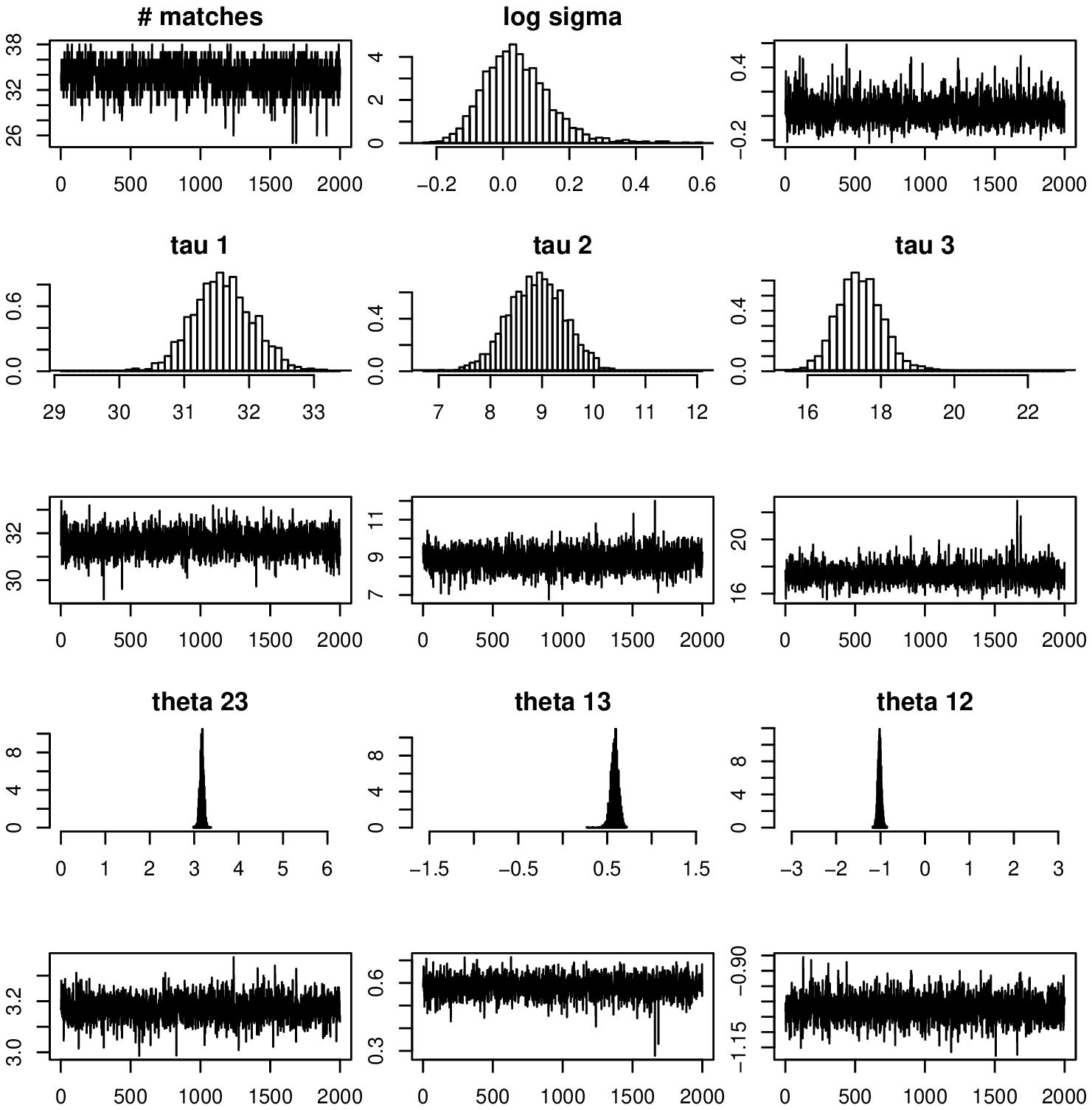}{Time series traces and histograms of the MCMC run
of Section \ref{sec:3deg}, based on a thinned sub-sample of 2000 after burn-in.}{5}

There is information about the identities of the amino acids 
in the two configurations: we defer use of this to Section \ref{sec:amino}.

The MCMC sampler described in Section \ref{sec:mcmc} was run for 
1 000 000 sweeps, after a burn-in period of 200 000 sweeps, considered on the
basis of an informal visual assessment of time series traces to be adequate
for convergence. 
Prior and hyperprior settings were: $\alpha=1$, $\beta=36$, $\mu_\tau=(0,0,0)^T$,
$\sigma_\tau=50.0$, $\lambda/\rho=0.003$ and the matrix $F_0$ defining the prior
for $A$ set to the zero matrix. The sampler parameter $p^\star$ was set to 0.5, and
we made updates to $M$ 10 times in each sweep.
Such a run took about 71 seconds on a 800MHz PC. Acceptance rates for the moves updating
$M$ were between 0.41\% and 5.6\%.

The posterior expectation and variance of $\tau$ were estimated to
be $(31.60,8.89,17.44)^T$ and
$$
\left(
\begin{array}{rrr}
 0.227 & 0.120 & -0.044 \\
 0.120 & 0.307 & 0.176 \\
 -0.044 & 0.176 & 0.428 \\
\end{array}
\right)
$$
The posterior mean and variance of $\sigma$ are 1.051 and 0.00996.
In representing the centre of the posterior
distribution for the rotation matrix $A$, we
we need to use a definition of mean appropriate to
the geometry. We form the mean elementwise
from a thinned sample of 2000 
values of $A$ from the post-burn-in MCMC run.
This mean matrix $\overline{A}$ is of course not a rotation matrix,
but post-multiplication by the positive definite symmetric
square root of $\overline{A}^T\overline{A}$
yields a rotation matrix that
is known as its polar part (see Mardia and Jupp, p. 286, 290).
This is an appropriate measure of location
of the posterior, and takes the value
\bel{meana1}
\widehat{A}=\left(
\begin{array}{rrr}
0.4339 & -0.8444 & 0.3140 \\
-0.7118 & -0.5350 & -0.4550 \\
0.5522 & -0.0261 & -0.8333 \\
\end{array}
\right)
\ee
in this case. 

The 40 most probable matches between $x$ and $y$ points are listed in
supplementary information;
there is no duplication in their indices until the 39th
match: $k=12$ also appears in the 38th match. 
We can conclude that for all values of $K$
greater than 0.2895 (the marginal posterior probability associated
with the 39th match), the optimal Bayesian matching is given
by an appropriate leading subset of the matches.
For example if this cost ratio is 0.5 we take the first 36
matches; these are displayed graphically
in Figure \ref{fig:nicolaconf42}; in this 3-dimensional example, the axes signify
the first two principle coordinates of the combined cloud of data.

As would be anticipated, simultaneous inference for the rotation $A$
and the matching matrix $M$ (as well as $\tau$ and $\sigma$)
is a considerably greater challenge
for MCMC than is the problem of the previous section, where 
the rotation matrix is held fixed. It is clear that there is a
possibility of severe multi-modality in the posterior, with the 
conditional posterior for $M$ and $\tau$ given $A$ depending 
strongly on $A$. This challenge is quantified empirically by a 
heavy-tailed distribution of times to convergence, and by `meta-stability'
in the time series plots of various monitoring statistics
against simulation time. We found the log-posterior
to be a useful summary statistic for quality of fit, and
pilot runs provided experience to choose a threshold value,
exceedance of which we hypothesised diagnosed convergence to
the main mode of the posterior.

To investigate
multimodality and convergence time, we conducted a study in which
the MCMC run described was repeated -- with the same parameters -- 
from 100 different initial configurations, obtained by independent
random rotations as initial settings for $A$. After short runs of
50 000 sweeps, we tested whether the threshold log-posterior value had
been exceeded, and if not the run was abandoned. 83 out of the 100
runs passed this test, and these were allowed to run on for a further 
450 000 sweeps. Every one of these 83 long runs provided exactly the
same set of 36 most probable matches, and we therefore felt justified to
conclude that they had not been trapped in a subsidiary mode of
the posterior, and that it was safe to draw inference from the
results. This conclusion is specific to the data set and
parameter settings used, and it would be straightforward to contrive
artificial data where multiple modes were more equal in probability
content. In such cases more sophisticated MCMC samplers would be needed.
 
\subsection{Prior settings and sensitivity}
\label{sec:sensitivity}

Our analysis depends of course on the settings of the
hyperparameters $\lambda/\rho$ (see Section \ref{sec:pp}), $F_0$ (Section \ref{sec:rotmat}),
and $\mu_\tau$, $\sigma_\tau$, $\alpha$, $\beta$ (Section \ref{sec:prior}).
These allow the provision of real prior information from the experimental context,
if it is available. 

For a default analysis in the absence of such information, we would
set $F_0$ to the zero matrix (a uniform prior on $A$),
$\mu_\tau$ to be the zero vector, and $\sigma_\tau$ of the order of twice the
distance between the centres of gravity of the two configurations.
We fix $\alpha=1$, giving an exponential prior distribution
for $\sigma^{-2}$.
Here we briefly discuss settings of, and sensitivity to, the remaining two
parameters, the scalars $\lambda/\rho$ and $\beta$.

Sensitivity to $\lambda/\rho$ is
pronounced, as might be anticipated. This parameter ratio
has a very direct role in determining whether an $(x_j,y_k)$ pair
are noisy observations of the same hidden $\mu_i$
point or not, after transformation, since it controls the density of
hidden points. In practice, we should not expect to
be able to draw inference about matching without
real prior knowledge about this ratio or an equivalent measure of
the prior tendency of points to be matched.

The prior for the number of matches $L$ is parameterised
by $\lambda/\rho$: see (\ref{lprior}).
This distribution is non-standard, but
very well-behaved. It is clear from inspection that
setting $\lambda/\rho$ equal to 
$(m-\Lbar)(n-\Lbar)/\Lbar v$ yields a mode of
$L$ that is within 1 of $\Lbar$, and numerical
calculation in the context of the example
in Section \ref{sec:3deg}, verifies that for all possible
`prior guesses' $\Lbar$ for $L$, the prior expectation and median
are also both equal to $\Lbar$ to the nearest integer.
Thus prior information about $L$ is directly informative
about the parameter ratio $\lambda/\rho$.
As long as $v$ is known, or at least a representative value
provided, and the analyst is able to make a prior
guess $\Lbar$ at the number of matches,
this suggests a reasonable way to specify $\lambda/\rho$.
The posterior distribution for $L$ tracks the prior
rather closely, confirming that the raw data carry little information
about the number of matches.

The hyperparameter $\beta$ is an inverse scale parameter for
the precision of the noise terms $\eps$; thus as $\beta$ increases,
we expect that $\sigma^2=\mbox{var}(\eps)$ increases too. The runs we have presented
used $\beta=36$; reducing this by a factor of 2 makes 
minimal difference to the posterior inference for 
either $\sigma^2$ or $M$. However, increasing $\beta$ by a factor of
2 leads to a 3-fold increase in $\sigma$ and a sharp reduction in the
number of matches -- the posterior expectation of $L$ goes down from
around 34 to 26. The latter observation is perhaps counter-intuitive,
until one realises that when $\sigma$ is larger, it becomes relatively
less likely that points that are nearly coincident (after transformation)
are in fact matched.

Finally, it would be desirable to assess the sensitivity to the Poisson assumption
for the hidden point model, but this would be extremely onerous to
do directly, since alternatives would require a substantially modified formulation
and implementation. There is scientific reason to doubt
the Poisson assumption; for
example, the minimum spacing between the centres of gravity of the amino acids
in proteins is approximately 3.8 Angstroms. However, 
experiments reported in Mardia, Nyirongo and Westhead (2005)
do at least suggest strongly that the ability of our method to
detect matches is little affected by real hard-core effects.


\subsection{Using information about types of amino acid}
\label{sec:amino}
The protein alignment data includes identifiers of the type of amino acid
at each point (see supplementary information). 
There are 20 different types, which can be categorised
into 4 groups: hydrophobic, charged, polar and glycine;
we use the group identifiers as colours 
in defining a modified likelihood as in Section \ref{modlik}.
The parameter values taken were
$\gamma=1.0$ and $\delta=-0.5$, providing a strong preference for
like-coloured matching ($\exp(\gamma-\delta)\approx 4.48$).
The analysis was repeated with this modified model, leaving all other
details unchanged.

The 40 most marginally probable matches are listed in supplementary information,
along with displayes of the optimal alignment.
The 36 most probable matches, which together form the optimal
matching whtn $K=0.4$,
are identical to those found in the previous section; however, there are modest variations in the posterior
probabilities attached to individual matches.

The posterior expectation and variance of $\tau$ were now estimated to
be $(31.94,8.94,17.61)^T$ (slightly shifted from that obtained in
the analysis of the previous section) and
$$
\left(
\begin{array}{rrr}
1.284 & -0.763 & -0.118 \\
-0.763 & 3.534 & -0.015 \\
 -0.118 & -0.015 &  1.320 \\
\end{array}
\right)
$$
The posterior mean and variance of $\sigma$ are 1.3122 and 0.1984.
The increased estimate of $\sigma$ is perhaps anticipated.
The centre of the posterior distribution of $A$
is in this case:
\bel{meana2}
\widehat{A}=\left(
\begin{array}{rrr}
0.4240 & -0.8512 & 0.3092 \\
-0.7235 & -0.5237 & -0.4497 \\
 0.5447 & -0.0331 & -0.8379 \\
\end{array}
\right).
\ee

\myfig{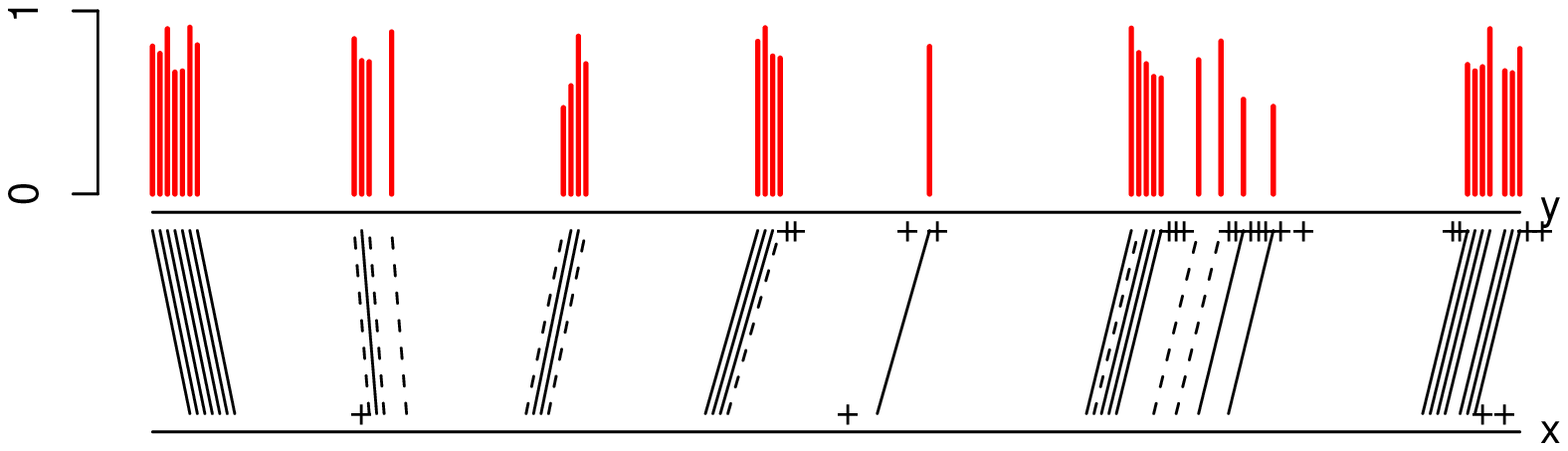}{The optimal matching (36 matches), when $K=0.4$, in the protein alignment analysis data,
using colouring information, with $\gamma=1.0$ and $\delta=-0.5$;
matches are signified by line segments joining the sequence number 
of the point in the $x$ configuration to that of the
matched point in the $y$ configuration.
The solid lines indicate the 27 matches identified by Gold et al (2002);
our method discovers all of these, together with the 9
further matches shown with broken lines.
The height of the vertical bars indicate the marginal
probabilities of each match.
The + symbols denote points that are present in either configuration 
but are not matched.}{6.5}

In the approach to the analysis of these data taken by Gold et al (2002),
the matching between the configurations 
was performed in two stages, and is not driven by an 
explicit probability model. First, inter-point distances
$d(\cdot,\cdot)$ were calculated within each configuration.
These distances are invariant under the rigid body motions
considered here. A maximal set of pairs of indices
$\{(j_1,k_1),(j_2,k_2),\ldots\}$, with no ties
among the $j$s or $k$s, is found such that 
$|d(x_{j_r},x_{j_s})-d(y_{k_r},y_{k_s})|$ is less than
some threshold, for all $s\neq r$. This is done using graph
theoretical algorithms of Bron and Kerbosch (1973) and 	
Carraghan and Pardos (1990), applied to a product graph
whose vertices are labelled with $(j,k)$ pairs.
This first stage of the matching alogrithm was formulated
by Kuhl et al (1984).

In the second stage, the matches are scored
using the amino acid information, assigning a score of 1
for identity of the amino acids, and 0.5 when the amino acids
are different but fall in the same group. The initial list
of matches from stage one is then permuted so as to
maximise the total score.

Once the matches are found the 	
rigid body transformation is estimated by Procrustes 
analysis; for example, see Dryden and Mardia (1998, pp 176-178).

It is interesting to compare the rotation
matrix resulting from this method, namely 
$$
A=\left(\begin{array}{rrr}
 0.441  & -0.841 & 0.312 \\
 -0.678 & -0.541 & -0.498 \\
  0.588 &  0.008 & -0.809 \\
\end{array}
\right)
$$
with that obtained by our method. The trace of the orthogonal
matrix taking $A$ to $\widehat{A}$ is approximately
$1+2\cos 0.07$, so the two differ by
a rotation of only 0.07 radians.
Figure \ref{fig:comb43}
provides a comparison between the matchings achieved by the two
approaches. Of the 27 matches identified by Gold et al, 
14 are among the most probable 20 that we find, and all 27 are 
among the first 35.


A referee has raised with us the role of sequence ordering along
the protein in inference about alignment and matching.
The example in this section concerns ligand binding site
matching, in which biologically relevant matches do not necessarily preserve
sequential ordering, in contrast to the more familiar 
situation of aligning protein backbones;
see for example Eidhammer et al (2004, pp. 333--334).
Examples are trypsin-subtilisin with
similar active sites and unrelated folds, and many adenine binding
sites in different folds. Somewhat remarkably, although sequence ordering
is not used in our analysis, the resulting matches do perfectly
respect this ordering. This is visualised in Figure \ref{fig:comb43},
which also reveals that some but not all of the matches revealed
by our analysis additional to those of Gold et al (2002) extend
already matched segments.
In this particular data set, the sites must come from very closely
related folds and would probably also be alignable by sequence-preserving
methods aligning full structures. Intriguingly, in this example at least,
knowledge of the sequence ordering would provide no additional
information beyond that extracted from the point coordinates and
amino-acid groups by our approach.

\section{Discussion}

The main conclusion of this paper is that a probability
model based approach is successful in allowing
simultaneous inference about
partial matching between two point configurations, and 
a geometrical transformation between the coordinate systems
in which the configurations are measured. This seems
an advance over previous more ad-hoc methods.

We have only used the translation and rigid motion groups in illustrating our
methodology. However, the formulation allows inference about various other group
transformations such as affine transformation, and so on. 
The fairly straightforward MCMC implementation presented here has
proved adequate for the models and data sets considered, although
allowing rotations did increase the needed run lengths considerably.
We anticipate that, at least for
models allowing rotations, dealing with larger data sets
will be much more challenging, since small 
rotational perturbations generate large displacements at sites far
from the axis of rotation; moves that simultaneously perturb
allocations and geometrical and error distribution parameters 
will be necessary for good performance. We also anticipate
more severe difficulties from multi-modality that were exposed
in Section \ref{sec:3deg}. 

An important task left for future work is a formulation that
allows smooth nonparametric transformations between coordinate
systems, setting warping into a model-based framework;
this would be important in dealing more comprehensively with
gel matching problems.

We have only used pairwise
comparisons but there is scope for taking multiple combinations such as
triads. The transformations considered above
are parametric but some non-parametric alternatives such as non-linear
deformations may be useful in some cases, e.g. to deal with dynamic aspects
of the atoms in a protein. We have
considered only two configurations but a natural extension
would be to take three or more point configurations
simultaneously, and make joint inference about patterns of matching between
the configurations and the various geometrical transformations involved.
More straightforward extensions would be
to allow for non-Gaussian noise, other types of prior and so on. 

Kent et al (2004) have treated the unlabelled case by using a 
different model. While matching two configurations, one of them 
is taken as the population and the second as a random sample 
from this population after an unknown transformation.
This approach is different from the symmetrical model for the 
two configurations proposed here.  Further the emphasis in 
Kent et al (2004) is on maximum likelihood inference using the EM-algorithm.

Recent independent work by Dryden, Hirst and Melville (2005),
addresses a similar problem of matching unlabelled point sets.
Their approach has some substantial differences, for example
there is assymmetry
in comparing two configurations, one being treated as a perturbation
of the other. The geometrical transformation parameters are given
uniform priors and maximised out, using standard ideas from
shape analysis, rather than integrated out as in our fully
Bayesian approach. Neither the loss function basis for
estimating matches, nor the treatment of partial labelling,
appear. 

There is other statistical work on alignment and matching in
proteins by Wu et al (1998) and Schmidler (2004),
which in contrast does use sequence information.
Further work is needed to clarify the
relationships between all these methods and their comparative
performance.

Finally, in the context of using methods such as ours in 
database search, often the reason for assessing protein alignment,
there are issues related to multiple comparisons. These are not discussed here,
but the answers will depend on the size of the database as well as the number of points in
the query site.

\section*{Acknowledgements}
We are grateful to Nicola Gold and Dave Westhead for their 
many helpful discussions, and in particular for the data in 
Example 2, and to Vysaul Nyirongo and Charles Taylor for 
various helpful comments. 

\section*{References}

\begin{list}{}{\setlength{\itemindent}{-\leftmargin}}

\item Abramowitz, M. and Stegun, I. A. (1970).
{Handbook of Mathematical Functions}. Dover, New York.

\item Bron, C. and Kerbosch, J. (1973).
Algorithm 457; finding all cliques of an undirected graph. 
{\it Communication of the ACM}, {\bf 16}, 575--577.

\item Carraghan, R. and Pardalos, P. M. (1990).
Exact algorithm for the minimal clique problem.  
{\it Operations Research Letters}, {\bf 9}, 375.

\item Chui, H. and Rangarajan, A. (2000). A new algorithm for non-rigid point matching. {\em IEEE Conference on Computer Vision and Pattern Recognition.} {\bf 2}, 44--51.

\item Cross, A. D. J. and Hancock, E. R. (1998). Graph matching with dual-step
{EM} algorithm.  {\em IEEE transactions on pattern analysis and machine
intelligence.} {\bf 20}, 1236--1253.

\item Downs, T. D. (1972). Orientation statistics.
{\it Biometrika}, {\bf 59}, 665--676.

\item Dryden, I. L., Hirst, J. D. and Melville, J. L. (2005).
Statistical analysis of unlabelled point sets: comparing
molecules in chemoinformatics.
Under revision for {\it Biometrics}.

\item Dryden, I. L. and Mardia, K. V. (1998). {\em Statistical shape analysis}.
Wiley, Chichester.

\item Eidhammer, T., Jonassen, T. and Taylor, W. R. (2004).
{\it Protein Bioinformatics}. 
Wiley, Chichester. 

\item Gold, N. D., Pickering, S. J., and Westhead, D. R. (2002). 
Protein functional site matching using graph theory techniques. 
In {\it Proceedings of
the International Conference on Bioinformatics}, Bangkok, Thailand,
page 79.


\item Green, P. J. (2001).
A Primer on Markov chain Monte Carlo, pp. 1--62 of 
{\it Complex Stochastic Systems}, Barndorff-Nielsen, O. E., Cox, D. R. and Kl\"{u}ppelberg, C. (eds.), Chapman and Hall, London.

\item Horgan, G. W., Creasey, A. and Fenton, B. (1992).
Superimposing two
dimensional gels to study genetic variation in malaria parasites. 
{\it Electrophoresis}, {\bf 13},871--875.

\item Kent, J. T., Mardia, K. V. and Taylor, C. C. (2004).
Matching problems for unlabelled configurations.
In {\it Bioinformatics, Images and Wavelets}, edited by
Aykroyd, R.G., Barber, S. and Mardia, K.V.
Proceedings of LASR 2004, Leeds University Press, Leeds.

\item Khatri, C. G. and Mardia, K. V. (1977). The von Mises--Fisher distribution
in orientation statistics.
{\it Journal of the Royal Statistical Society}, B, {\bf 39}, 95--106.

\item Kuhl, F. S., Crippen, G. M. and Friesen, D. K. (1984).
A combinatorial algorithm for calculating ligand binding.
{\it Journal of Computational Chemistry}, {\bf 5}, 24--34.

\item Mardia, K. V. and El-Atoum, S. A. M. (1976).
Bayesian inference for the von Mises--Fisher distribution. 
{\it Biometrika}, {\bf 63}, 203--205. 

\item Mardia, K. V. and Gadsden, R. J. (1977).
A circle of best fit for spherical data and areas of  
vulcanism. 
{\it Applied Statistics},{ \bf 26}, 238--245.

\item Mardia, K. V. and Jupp, P. E. (2000).
{\it Directional Statistics}, Wiley, Chichester.

\item Mardia K. V., Taylor, C. C, and Westhead, D. R. (2003). Structural bioinformatics revisited. In {\it LASR2003}, pp11--18. Leeds University Press.

\item Mardia, K. V., Nyirongo, V., and Westhead, D.R. (2005).
EM algorithm, Bayesian and distance approaches to matching active sites
{\it Mathematical and Statistical Annual Meeting in Bioinformatics}, Rothamsted, March 2005, Abstracts pp13-14.

\item Pedersen, L. (2002). {\em Analysis of two-dimensional electrophoresis gel images.} Ph.D thesis, IMM Technical University of Denmark.

\item Raffenetti, R. C. and Ruedenberg, K. (1970).
Parameterization of an orthogonal matrix
in terms of generalized Eulerian angles.
{\it International Journal of Quantum Chemistry}, {\bf IIIS}, 625--634.

\item Richardson, S. and Green, P. J. (1997). 
On Bayesian analysis of mixtures with an unknown number of components 
(with discussion). 
{\it Journal of the Royal Statistical Society}, B, {\bf 59}, 731--792.

\item Schmidler, S. C. (2004).
{\it Bayesian shape matching and structural alignment}.
Presentation at the 6th World Congress of the Bernoulli Society,
Barcelona, July 2004.

\item Wu, T. D., Schmidler, S. C., Hastie, T. and Brutlag, G. (1998).
Regression analysis of multiple protein structures. 
{\it Journal of Computational Biology}, {\bf 5}, pp 585--595.

\end{list}

\end{document}